\documentclass[12pt,a4paper]{article}
\usepackage{amsmath, amssymb, theorem, latexsym, epsfig}
 \usepackage{graphics}

\newtheorem{theorem}{Theorem}[section]
\newtheorem{lemma}[theorem]{Lemma}
\newtheorem{proposition}[theorem]{Proposition}
\newtheorem{definition}[theorem]{Definition}

\newtheorem{conjecture}[theorem]{Conjecture}
\newtheorem{corollary}[theorem]{Corollary}

\allowdisplaybreaks

\def\neweq#1{\begin{equation}\label{#1}}
\def\endeq{\end{equation}}
\numberwithin{equation}{section}
\begin{document}
\setlength{\parindent}{0 cm}
{\centering
\bfseries
{\LARGE Construction of Permutation Snarks}

\bigskip
\mdseries
{\large Jonas H\"agglund and Arthur Hoffmann-Ostenhof}\footnote{supported by the FWF project P20543.}

}

\begin{abstract}

A permutation snark is a snark which has a $2$-factor $F_2$ consisting of two chordless circuits; $F_2$ is called the permutation $2$-factor of $G$. We construct an infinite family $\mathcal H$ of cyclically $5$-edge connected permutation snarks. Moreover, we prove for every member $G \in \mathcal H$ that the permutation $2$-factor given by the construction of $G$ is not contained in any circuit double cover of $G$. 

\end{abstract}


\vspace{0.5cm}

\section{Introduction and main result}

A \emph{circuit} is defined to be a $2$-regular $2$-connected graph. A \emph{circuit double cover} (CDC) of a cubic graph $G$ is a set $S$ of circuits of $G$ such that every edge of $G$ is covered by exactly two circuits of $S$. A $2$-regular subgraph $D$ of $G$ is said to be contained in $S$ if every circuit of $D$ is an element of $S$.

A cubic graph $G$ with a $2$-factor $F_2$ which consists of two chordless circuits is called a \emph{cycle permutation graph} and
$F_2$ is called the \emph{permutation $2$-factor} of $G$. If $G$ is also a snark, then we say $G$ is a \emph{permutation snark}.
The Petersen graph has been for a long time the only known cyclically $5$-edge connected permutation snark.
In \cite{Br} twelve new cyclically $5$-edge connected permutation snarks have been discovered by computer search. Here, we present the first infinite family of cyclically $5$-edge connected permutation snarks.

We state the main theorem, see Theorem \ref{!!} and Corollary \ref{wichtig}.

\begin{theorem}
 For every $n \in \mathbb N$, there is a cyclically $5$-edge connected permutation snarks $G$ of order $10+24n$. Moreover, $G$ has a permutation $2$-factor which is not contained in any CDC of $G$.
\end{theorem}

Applying the above theorem we obtain infinitely many counterexamples to the following conjectures.

\begin{conjecture}\label{c2}\cite{Z1}
Let $G$ be a cyclically $5$-edge-connected cycle permutation graph. If $G$ is a snark, then $G$ must be the Petersen graph.
\end{conjecture}

\begin{conjecture}\label{c1}\cite{F}
If $G$ is an essentially $6$-edge-connected $4$-regular
graph with a transition system $T$, then $(G, T)$ has no compatible cycle 
decomposition if and only if $(G, T)$ is the bad loop or the bad $K_5$.
\end{conjecture}


\begin{conjecture}\label{c4} \cite{BJ, Z1} 
Let $G$ be a cyclically $5$-edge-connected cubic graph and $D$ be a set of pairwise disjoint circuits of $G$. Then $D$ is a
subset of a CDC, unless $G$ is the Petersen graph.
\end{conjecture}

\begin{conjecture}\label{c5} \cite{Z1} 
Let $G$ be a cycle permutation graph with the cordless circuits $C_1$ and $C_2$ where $C_1 \cup C_2$ is a 2-factor. If $G$ is cyclically $5$-edge-connected
and there is no CDC which includes both $C_1$ and $C_2$, then $G$ must be the Petersen graph.
\end{conjecture}

Note that finitely many counterexamples to the above conjectures were found in \cite{Br} by computer search.

\section{Definitions and proofs}
We refer to \cite {Z1} for the definition of a multitpole and an half-edge.
Moreover, for terminology not defined here we refer to \cite{Bo}. We use a more general definition of a CDC than the one stated in the introduction.

\begin{definition}\label{pcdc}
We say a set $S=\{A_1,A_2,...,A_m\}$ is a path circuit double cover (PCDC) of a graph $G$ if the the following is true

1. $A_i$ is a subgraph of $G$ where every component of $A_i$ is either a circuit or a path with both endvertices being vertices of degree $1$ in $G$, $ \forall \,i \in \{1,2,...,m\}$.

2. $\sum_{i=1}^m |e \cap E(A_i)|=2 \,\,\, \forall \,e \in E(G)$. 

\end{definition}

If no $A_i$ contains a path as a component, then we call $S$ a \emph{CDC} of $G$ and if $|S|=k$, then we call $S$ a \emph{$k$-CDC} of $G$. Obviously, a PCDC is a CDC if $G$ contains no vertex of degree $1$. For a survey on CDC's, see \cite{Z1,Z2}.

Later we need the following known lemma \cite{Z1}.

\begin{lemma}\label{bla}
Let $G$ be a $3$-edge colorable cubic graph and $D$ be a $2$-regular subgraph of $G$. Then $G$ has a $4$-CDC $S$ with $D \in S$.  
\end{lemma}

\begin{definition}
Let $A \in S $ be given where $S$ is a PCDC of a graph $G$. Let $e$ be an half-edge or edge of $A$, then $[e]$ denotes the 
unique element of $S$ which contains $e$ and which is not $A$. We say $[\,]$ refers to $A$.
\end{definition}

\begin{definition}\label{c5}
Let $Q^i$ with $i \in \{1,2,...,4\}$ be a cyclically $5$-edge connected permutation snark with a permutation $2$-factor $F^i$ such that $F^i$ is not contained in any CDC of $Q^i$. The two circuits of $F^i$ are denoted by $C^i_1$ and $C^i_2$. We may assume w.l.o.g. that $C^i_1$ ($C^i_2$) contains a subpath which has the following vertices in the following consecutive order: $x^i_1$, $x^i_2$, $z^i_2$, $x^i_6$ ($x^i_4$, $z^i_1$, $x^i_5$), such that $z^i_1z^i_2 \in E(Q^i)$.
\end{definition}

\begin{definition}
Let $\widetilde{ Q }^{i}$ be the graph which is obtained from $Q^i$ by removing the edge $x^i_1x^i_2$, the vertices 
$z^i_1$, $z^i_2$ and by adding the vertex $y^i_j$, $j=1,2,...,6$ and the edges 
$e^i_3:=x^i_2y^i_3$, $e^i_s:=x^i_sy^i_s$, $s=1,2,4,5,6$, see Figure 1. 
\end{definition}

By an \emph{end-edge} of a graph $G$, we mean an edge which is incident with a vertex of $G$ with degree $1$ in $G$. 

\begin{definition}\label{A}
The six end-edges of $\widetilde{Q}^i$ together with the remaining edges of $F^i$ in $Q^i$ induce the following three paths in $\widetilde{Q}^i$:
the path $A^i_1$ with end-edges $e^i_1$ and $e^i_6$,
the path $A^i_2$ with end-edges $e^i_4$ and $e^i_5$ and the path $A^i_3$ with end-edges $e^i_2$ and $e^i_3$.
Moreover, set $A^i:= A^i_1 \cup A^i_2 \cup A^i_3$ and $\mathbb A^i := \{A^i_1, A^i_2, A^i_3 \}$.
 \end{definition}

\begin{figure}[htpb]
\centering\epsfig{file=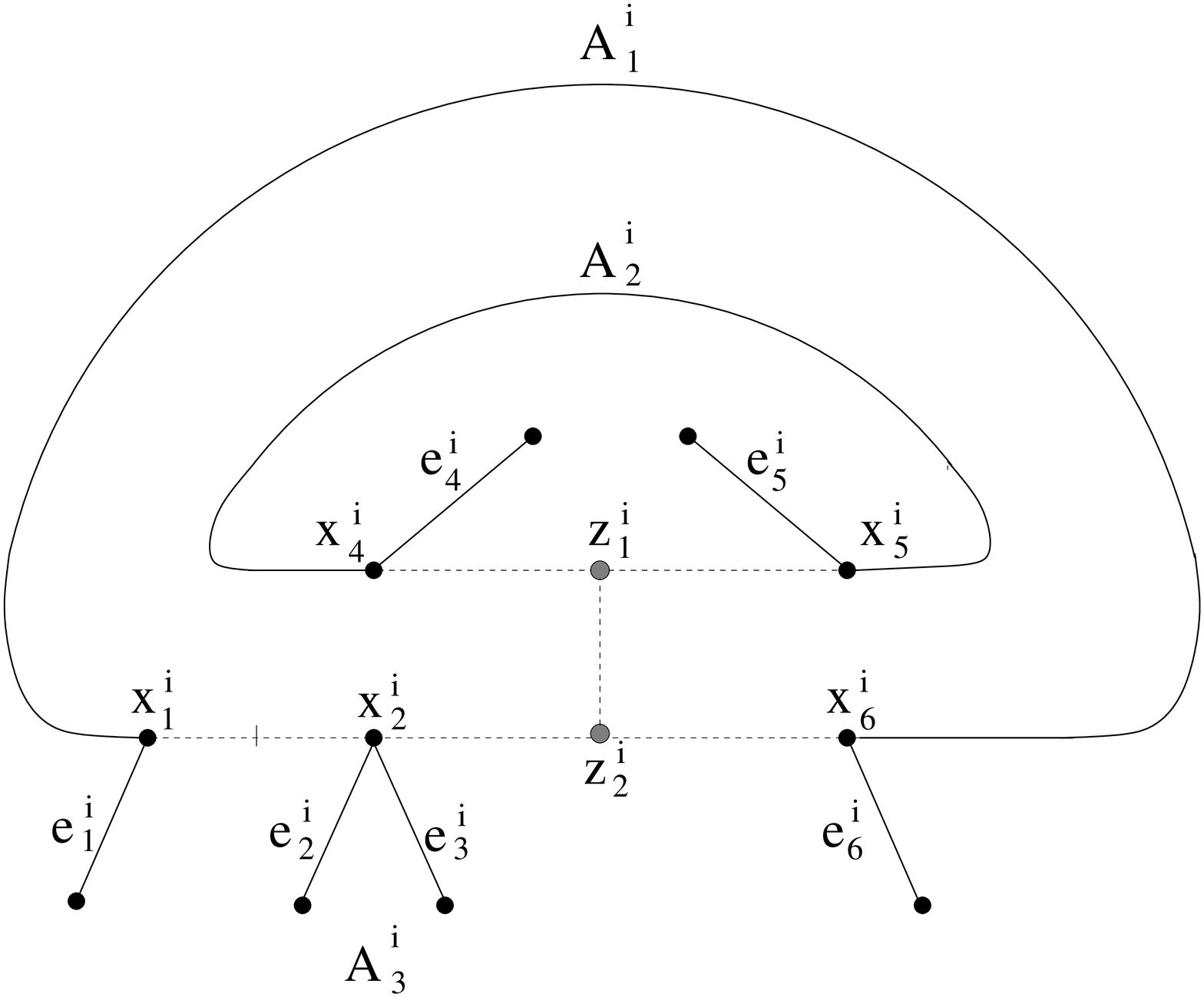,width=3.2in}
\caption{The graph $\widetilde{Q}^i$, $i \in \{1,2,3,4\}$.} \label{color1}
\end{figure}

We recall that $[\,]$ refers to the given element of a PCDC or CDC. We need the following propositions and lemmas for proving Theorem \ref{!!}.

\begin{proposition}\label{qa}
Let $\widetilde{Q}^i$ and $A^i$ with $i \in \{1,2,...,4\}$ be defined as above.
Then every PCDC $S$ of $\widetilde{Q}^i$ with $A^i \in S$ satisfies the following.

$(1)\,\,$ If $[e^i_4] \not= [e^i_5]$, then $[e^i_1] \not\in \{[e^i_2],[e^i_3]\}$. If $[e^i_1] \in \{[e^i_2],[e^i_3]\}$, then $[e^i_4] = [e^i_5]$.\\
$(2)\,\,$ $[e^i_2]\not=[e^i_3]$.
\end{proposition}

Proof. Suppose by contradiction that one of the two conclusions of $(1)$ is not fulfilled. Then $S$ implies a CDC of $Q^i$ containing $F^i$ which contradicts the definition of $Q^i$. If $(2)$ is not fulfilled, then the edge $e \not\in A^i_3$ which is incident with $x^i_2$ cannot be covered 
by $S$ which is impossible. Hence, the proof is finished.

\begin{corollary}\label{ww}
Let $S$ with $A^i \in S$ be a PCDC of $\widetilde{Q}^i$ with $i \in \{1,2,...,4\}$. If $[e^i_4] \not= [e^i_5]$, then

$(1)\,\,$ $|\{[e^i_1],[e^i_2],[e^i_3]\}|=|\{[e^i_4],[e^i_5],[e^i_6]\}|=3$.\\
$(2)\,\,$ $|\{[e^i_1],[e^i_2],...,[e^i_6]\}|=3$.
\end{corollary}

Proof. Since by Proposition \ref{qa}, $[e^i_2] \not= [e^i_3]$ and $[e^i_1] \not\in \{[e^i_2],[e^i_3]\}$, the corollary follows.

\begin{definition}
Denote by $P^i$ with $i \in \{1,2,3,4\}$ the connected multipole which is obtained from $\widetilde{Q}^i$ by transforming every end-edge of $\widetilde{Q}^i$ 
into an half-edge except for $i=1$, $e^1_2$; for $i=2$, $e^2_3$; for $i=3$, $e^3_2$; and for $i=4$, $e^4_3$.  
\end{definition}

\begin{definition}\label{Hcon}
Denote by $H(Q^1,Q^2,Q^3,Q^4)$ or in short by $H$ the cubic graph which is constructed from $P^1$, $P^2$, $P^ 3$ and $P^4$, as illustrated in Figure 2, by gluing together half-edges, identifying vertices of degree $1$ and by adding the edge $\alpha$.  
\end{definition}

Note that we keep in $H$ the edge labels of $P^i$, respectively, of $Q^i$, see Figure 2.

\begin{definition}
Let $A \in \mathbb A^i$ (Def. \ref{A}), i.e. $A \subseteq \widetilde {Q}^i$, $i \in \{1,2,3,4\}$. Then $A \subseteq H$ 
is defined to be the path in $H$ containing all edges of $H$ which have the same edge-labels as $A \subseteq \widetilde {Q}^i$; if an edge $e$ of 
$A \subseteq \widetilde {Q}^i$ corresponds to an half-edge of $H$ then the edge of $H$ which contains $e$ is defined to be part of $A \subseteq H$.
\end{definition}

\begin{definition}\label{c1c2}
Denote by $F$ the permutation $2$-factor of $H$ with $F:=A^1 \cup A^2 \cup A^3 \cup A^4$ and $A_i \subseteq H$, $i=1,2,3,4$, see Figure 1, 2 and 3.
\end{definition}

\begin{figure}[htpb]
\centering\epsfig{file=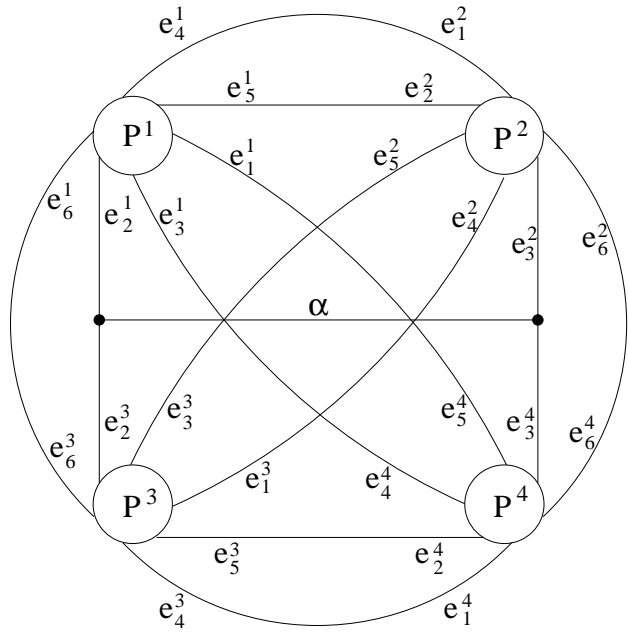,width=3.2in}
\caption{The graph $H$.} \label{color1}
\end{figure}

\begin{figure}[htpb]
\centering\epsfig{file=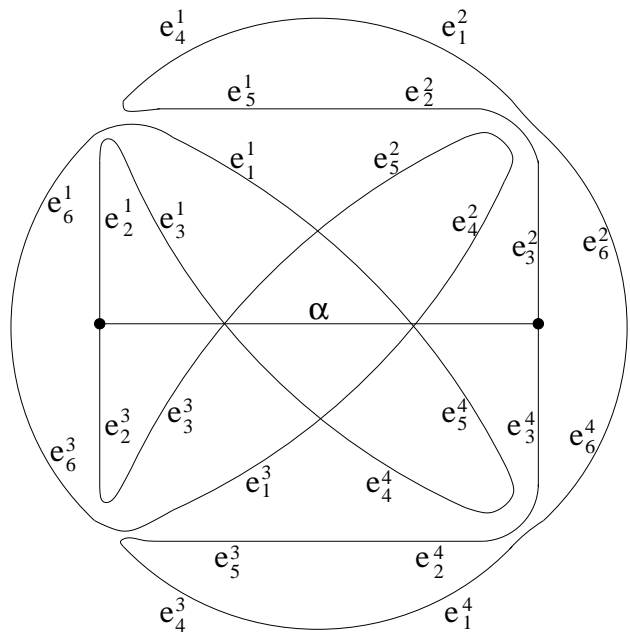,width=3.2in}
\caption{The permutation $2$-factor $F$ of $H$ and $\alpha \not\in E(F)$.} \label{ff}
\end{figure}

\begin{lemma} \label{1}
Let $S$ be a CDC of $H$ with $F \in S$. Then $[e^i_4] \not= [e^i_5]$ for some $i \in \{1,2,3,4\}$. 
\end{lemma}

Proof by contradiction. Consider $P^3$ in Figure 2. Since $[e^3_4]=[e^3_5]$ we obtain $[e^4_2]=[e^4_1]$. Therefore and since $[e^4_4]= [e^4_5]$
it follows that $[e^4_3]= [e^4_6]$. Consider $P^1$. By analogous arguments, $[e^2_3]= [e^2_6]$. Since $[e^2_6]= [e^4_6]$ we obtain $[e^4_3]= [e^2_3]$ which is impossible.

\begin{lemma} \label{2} 
Let $S$ be a CDC of $H$ with $F \in S$ and let $[e^i_4] \not= [e^i_5]$ for some $i \in \{1,2,3,4\}$. Then $[e^i_4] \not= [e^i_5]$ $\forall i \in \{1,2,3,4\}$. 
\end{lemma}

Proof. Let $[e^1_4] \not= [e^1_5]$. Then by Proposition \ref{qa} $(1)$, $[e^1_1] \not= [e^1_3]$. Hence $[e^4_4] \not= [e^4_5]$ and by 
Proposition \ref{qa} $(1)$, $[e^4_1] \not= [e^4_2]$. Hence $[e^3_4] \not= [e^3_5]$ and thus by Proposition \ref{qa} $(1)$, $[e^3_1] \not= [e^3_3]$ implying $[e^2_4] \not= [e^2_5]$. Each of the three remaining cases to consider, i.e. $[e^i_4] \not= [e^i_5]$, $i=2,3,4$, can be proven analogously.

Corollary \ref{ww}, Lemma \ref {1} and Lemma \ref{2} imply the following proposition.

\begin{proposition}\label{!}
Let $S$ be a CDC of $H$ with $F \in S$. Then the following is true for $i=1,2,3,4$.

$(1)\,\,$ $[e^i_4] \not= [e^i_5]$.\\
$(2)\,\,$ $|\{[e^i_1],[e^i_2],[e^i_3]\}|=|\{[e^i_4],[e^i_5],[e^i_6]\}|=3$.\\
$(3)\,\,$ $|\{[e^i_1],[e^i_2],...,[e^i_6]\}|=3$.
\end{proposition}

For the proof of the next theorem we form a new cubic graph $H'$ from $H$. Consider for this purpose $P^i$, $i=1,2,3,4$ in Figure 2 as a vertex of degree $6$ and split every $P^i$ into two vertices $v^i$ and $w^i$ such that $v^i$ ($w^i$) is incident with $e^i_j$, $j=1,2,3$ ($\,e^i_j$, $j=4,5,6$) to obtain $H'$. For reasons of convenience we do not use the edge-labels of $H$ for $H'$, see Figure 4.

\begin{figure}[htpb]
\centering\epsfig{file=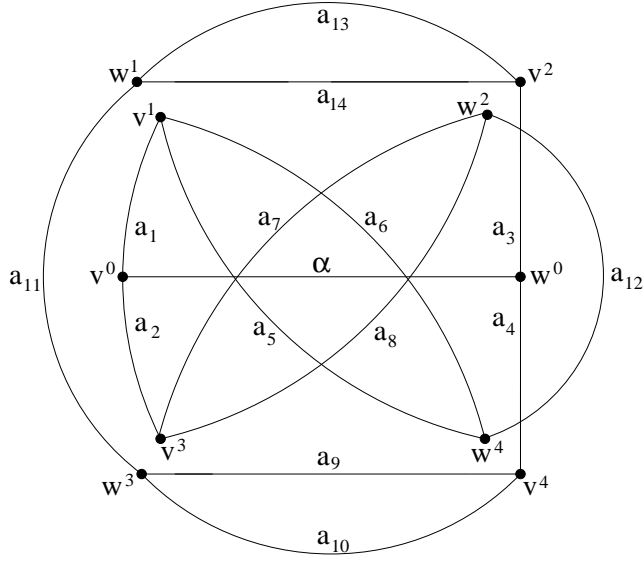,width=3.3in}
\caption{The graph $H'$.} \label{color1}
\end{figure}

\begin{theorem}\label{!!}
Let $S$ be CDC of $H$, then $F \not\in S$. 
\end{theorem}

Proof. Note that for every CDC $S$ of $H$ with $F \in S$, $\{[e^1_2],[e^3_2]\}=\{[e^2_3],[e^4_3]\}$.

$(2)$ and $(3)$ in Proposition \ref{!} imply that it suffices to show that there is no proper edge-coloring $f: \{E(H')- \alpha\} \mapsto \mathbb N$ such that $v^i$ and $w^i$ in $H'$ are incident with the same colors, $\forall \,i \in \{0,1,2,3,4\}$. Let $v \in V(H')$, then $E_v$ denotes the edge-set containing all edges of $H'$ incident with $v$. We proceed by contradiction. There are two cases to consider.

Case 1. $f(a_1)= f(a_3)=1$ and $f(a_2)= f(a_4)=2$.\\
Since $f(a_1)=1$ and $a_1 \in E_{v^1}$ there is $x \in E_{w^1}$ with $f(x)=1$. Since $f(a_3)=1$, $f(a_{11})=1$. Since $f(a_2)=2$ and since $a_2 \in E_{v^3}$ there is $y \in E_{w^3}$ with $f(y)=2$ which is impossible since $f(a_{11})=1$ and $f(a_4)=2$.

Case 2. $f(a_1)= f(a_4)=1$ and $f(a_2)= f(a_3)=2$.\\
Since $f(a_3)=2$ and $a_3 \in E_{v^2}$ there is $x \in E_{w^2}$ with $f(x)=2$. Since $f(a_2)=2$, $f(a_{12})=2$. 
Since $f(a_4)=1$ and $a_4 \in E_{v^4}$ there is $y \in E_{w^4}$ with $f(y)=1$ which is impossible since $f(a_1)=1$ and $f(a_{12})=2$.

The \emph{cyclic edge-connectivity} of a graph $G$ which contains two vertex-disjoint circuits is denoted by $\lambda_c(G)$; it is the minimum number of edges one needs to delete from $G$ in order to obtain two components such that each of them contains a circuit. In order to show that $\lambda_c (H) > 4$, we need several results.

\begin{definition}
Let $G$ be a graph with a given $2$-factor $F_2$ consisting of two chordless circuits. We call $e \in E(G)$ a spoke if $e \not\in E(F_2)$.
\end{definition}

\begin{proposition}\label{perch}
Let $G$ be a cubic graph with a $2$-factor $F_2$ consisting of two chordless circuits $C_1$, $C_2$.

(1) Let $|V(G)| \geq 8$, then $\lambda_c (G) \geq 4$.\\
(2) Let $|V(G)| \geq 10$. Then every cyclic $4$-edge cut $E_0$ of $G$ is a matching and $|E_0 \cap C_i|=2$, $i=1,2$.   
\end{proposition}

Proof. Suppose by contradiction that $M$ is a cyclic $3$-edge cut of $G$. Obviously, $M$ is matching. First we show that $M$ contains no spoke and consider two cases.

\emph{Case 1.} $M$ contains two or three spokes. \\ It is straightforward to see that $G-M$ is connected and thus this is impossible.

\emph{Case 2.} $M$ contains exactly one spoke. \\If $|M \cap E(C_i)| =1$ for $i=1,2$ then $C_1-M \subseteq G$ is a path which is connected by more than one spoke to $C_2-M \subseteq G$. Thus $G-M$ is connected. Hence w.l.o.g. $|M \cap E(C_1)|=2$. Since both paths of $C_1-M$ contain more than one vertex, both paths in $G$ are connected by more than one spoke to $C_2$ and thus $G-M$ is connected. Hence $M$ contains no spoke. 

Suppose $|M \cap E(C_1)|=2$ and thus $|M \cap E(C_2)|=1$. Then both paths of $C_1-M \subseteq G$ are connected by more than one spoke to the path $C_2-M \subseteq G$. Hence $G-M$ is connected which contradicts the assumption and thus finishes the first part of the proof.

Suppose $E_0$ contains a spoke $s$. For every subdivision $Z'$ of a cubic graph $Z$, $\lambda_c (Z')= \lambda_c (Z)$. Thus and by the first statement of the Proposition, $G-s$ is cyclically 4-edge connected. Hence $s \not\in E_0$. \\
Suppose $E_0$ is not a matching and let $a_1$ be adjacent with $a_2$ where $\{a_1,a_2\} \subseteq E_0 \cap E(C_1)$. Let $s$ be the unique spoke which is adjacent with $a_1$ and $a_2$. Then $E'_0:= E_0 - a_1 \cup s$ is a cyclic $4$-edge cut of $G$. This contradicts the previous observation that a spoke is not contained in any cyclic $4$-edge cut. Hence $E_0$ is a matching. Suppose $E_0$ contains no (one) edge of $C_1$ and thus four (three) edges of $C_2$. 
$C_2-E_0$ consists of four (three) paths where each of them is connected by a spoke to $C_1- E_0$ which is connected in both case. Hence $G-E_0$ is connected which contradicts the assumption and thus finishes the proof.

Let $V'$ be a subset of vertices of a graph $G$, then we denote by $\langle V'\rangle $ the \emph{vertex induced subgraph} of $G$.

\begin{lemma} \label{4cc}  
Let $G$ with $|V(G)| \geq 10$ be a cubic graph which contains a $2$-factor $F_2$ consisting of two chordless circuits $C_1$,$C_2$. Then the following is true and the analogous holds for $C_2$.

(1) $\lambda_c(G)= 4$ if and only if $C_1$ contains a path $L_1$ with $1 < |V(L_1)| < |V(C_1)|-1 $ such that $L_2:=\,\langle N(V(L_1)) \cap V(C_2) \rangle  $ is a 
path of $C_2$. In particular, the four distinct end-edges of $F_2-E(L_1)-E(L_2)$ form a cyclic $4$-edge cut of $G$.

(2) Let $E_0:= \{a_1,a_2,b_1,b_2\}$ be a cyclic $4$-edge cut of $G$ where by Prop. \ref{perch}, w.l.o.g. $\{a_1,a_2\} \subseteq E(C_1)$ and  $\{b_1,b_2\} \subseteq E(C_2)$. Denote the two paths of $C_1-a_1-a_2$ ($C_2-b_1-b_2$) by $L'_1$ and $L''_1$ ($L'_2$ and $L''_2$). Then, $\{\,\langle N(V(L'_1)) \cap V(C_2)\rangle , \, \langle N(V(L''_1)) \cap V(C_2) \rangle \,\} = \{L'_2, L''_2\}$.
\end{lemma}

Proof. We first prove (2). The four paths $L'_1, L''_1, L'_2, L''_2$ decompose the graph $C_1 \cup C_2 - E_0$. Every neighbor of a vertex of $L'_1$ or $L''_1$ in $C_2$ is thus contained in $L'_2$ or $L''_2$. Hence the equality in (2) does not hold if and only if $L'_1$ or $L''_1$ contains two distinct vertices $x, y$ such that $N(x) \cap V(C_2) \in V(L'_2)$ and $N(y) \cap V(C_2) \in V(L''_2)$. Let w.l.o.g.
$\{x,y\} \subseteq V(L'_1)$ and suppose by contradiction that $x$ and $y$ have the before described properties. Then $L'_1 \subseteq G-E_0$ is connected to $L'_2$ and $L''_2$. Since $L''_1 \subseteq G-E_0$ is connected to $L'_2$ or $L''_2$, $G-E_0$ is connected which contradicts the assumption and thus finishes this part of the proof.
 
We prove (1). Let $\lambda_c(G)= 4$ and let $E_0$ be a cyclic $4$-edge cut of $G$ as defined in (2). By Proposition \ref{perch}, $\{a_1,a_2\}$ is a matching of $G$. Hence the inequality in (1) is satisfied by setting $L_1:= L_1'$. It is straightforward to check that the four end-edges of 
$F_2-E(L_1)-E(L_2)$ form a cyclic $4$-edge cut of $G$. Hence the proof is finished.

\begin{lemma} \label{kkj} 
Let $A \subseteq H $ and $A \in \mathbb A^i $, $i \in \{1,2,3,4\}$. Then there is no cyclic $4$-edge cut $E_0 $ of $H$ such that $|E_0 \cap E(A)| =2$. 
\end{lemma}

Proof by contradiction. Set $E_0 \cap E(A)= \{a_1,a_2\}$. By Proposition \ref{perch}, $E_0$ is a matching. Hence $A \not= A^i_3$. There are two cases to consider.

\emph{Case 1.} $A= A^i_1$. Set $B:= A^i_2$. $A$ and $B$ belong to different components of $F \subseteq H$, see Figure 1 and Figure 3. 
Denote by $A^*$ the unique path of $A-a_1-a_2$ which connects one endvertex of $a_1$ with one endvertex of $a_2$.
Since $E_0$ is a cyclic $4$-edge cut and by Lemma \ref{4cc} (2) and by the structure of $P^i$, $B^*:=\langle N(V(A^*)) \cap V(B) \rangle $ is a subpath of $B$, see Figure 1. Since $\{a_1,a_2\}$ is a matching of $H$, $|V(A^*)| > 1$. Since $|V(A^*)| = |V(B^*)|$, $|V(B^*)| > 1$.

Denote by $w$ the neighbor of $x^i_2$ in $B$. Consider the graph $Q^i$ which was defined for constructing $P^i$, see Figure 1. Then $B^*$ is a path of $C^i_2 \subseteq Q^i$ (Def. \ref{c5}) with $1< |V(B^*)|< |V(C^i_2)|-1$ since $\{z^i_1, w\} \cap V(B^*) = \emptyset $. Moreover, $A^* \subseteq Q^i$ is a path of $C^i_1$ and every vertex of $A^*$ is adjacent to a vertex of $B^* \subseteq Q^i$ and vice versa. Therefore and by Lemma \ref{4cc} (1), $\lambda_c (Q^i)=4$ which is a contradiction to Def. \ref{c5}.

\emph{Case 2.} $A= A^i_2$. Set $B:= A^i_1$. Let $A^*$ be the unique path of $A-a_1-a_2$ which connects one endvertex of $a_1$ with one endvertex of $a_2$. Let $C$ denote the component of $F \subseteq H$ such that $A^* \not\subseteq C$. Since $E_0$ is a cyclic $4$-edge cut and by Lemma \ref{4cc} (2), $B^*:=\langle N(V(A^*)) \cap V(C) \rangle $ is a subpath of $C$. Denote by $w$ the neighbor of $x^i_2$ in $A$. Suppose $w \in V(A^*)$. Then $x^i_2 \in V(B^*)$. Since $|V(B^*)| > 1$, $B^*$ contains a vertex $y \not= x^i_2$. Since $N(V(A^*)-w) \cap V(C) \subseteq V(B)$ and since $x^i_2$ is not adjacent to a vertex of $B$, $B^*$ is not a path which is a contradiction. Hence, $w \not\in V(A^*)$. Similar to Case 1, $A^* \subseteq Q^i$ is a path of $C^i_2$ with  $\{z^i_1, w\} \cap V(A^*) = \emptyset$ where every vertex of $A^*$ is adjacent to a vertex of the path $B^* \subseteq Q^i$ and vice versa. Hence, by applying the same arguments as in Case 1, the proof is finished.

We need the following observations and notations for the proof of Theorem \ref{pppp}.

\begin{definition}
Let $X$ be a path with $|V(X)|>2$. Then $oX$ is the subpath of $X$ which contains all vertices of $X$ except the two endvertices of $X$. 
\end{definition}

Let $R$ be a set of subgraphs of $H$, then $E(R)$ denotes the union of the edge-sets of the subgraphs of $R$. 
Note that the vertices of the graph $J$ defined below are not the vertices of $H'$.

\begin{figure}[htpb]
\centering\epsfig{file=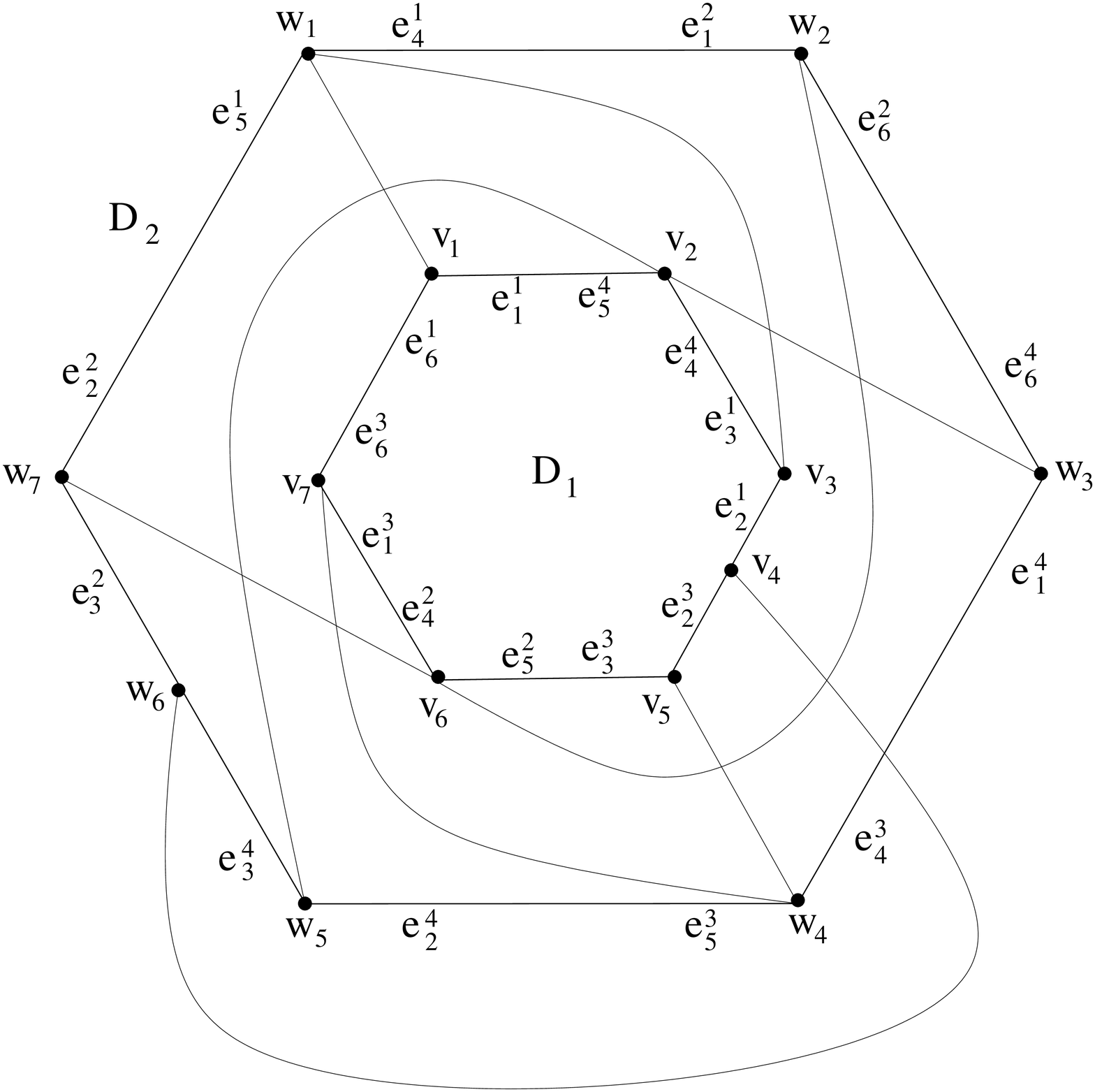,width=3.6in}
\caption{The graph $J$.} \label{inp}
\end{figure}

\begin{definition}\label{JJ}
Let $C_1$, $C_2$ denote the two circuits of $F \subseteq H$ where $e^1_2 \in E(C_1)$, see Figure 3.
Set $R:= \{oA^i_j|\,i \in \{1,2,3,4\} \textrm { and } j \in \{ 1,2,3 \} \}$. Let $J:=H/R$ be the the graph which is obtained from $H$ by contracting every element of $R$ in $H$ to a distinct vertex and by then replacing every multiple-edge by a single edge, see Figure \ref{inp}. The edges of $E(F)-E(R)$ induce a $2$-factor of $J$ which consists of two circuits $D_1$, $D_2$, see Figure 5.  
\end{definition}
  
Note that only subpaths of $F \subseteq H$ are contracted in the transformation from $H$ into $J$ and that $C_k \subseteq H$ is transformed into $D_k \subseteq J$, $k=1,2$. We keep the labels of the edges respectively half-edges of $E(F)-E(R) \subseteq E(H)$ for $E(D_1) \cup E(D_2) \subseteq E(J)$, see Figure 5. Every $oA^i_j \subseteq H$ corresponds to a vertex of $J$ and $v_4w_6 \in E(J)$ corresponds to $\alpha \in E(H)$. 
Set $\mathbb H:= \{V(oA^i_j) \,|\, i=1,2,3,4 \textrm { and } j=1,2,3\} \cup \{v_4\} \cup \{w_6\}$. Then $\mathbb H$ is a vertex partition of 
$V(H)$. Every $v \in V(J)$ corresponds to an element $\hat v$ of $\mathbb H$ and vice versa.

Let $h: V(H) \mapsto V(J)$ be the mapping where $h(x)$ is defined to be the unique $v \in V(J)$ such that $x \in \hat v$.\\  
Let $e \in E(C_k)-E(R)$, then $h'(e)$ denotes the corresponding edge of $D_k \subseteq J$, $k =1,2$.\\
Let $X$ be a subpath of $C_k \subseteq H$, then $X/R$ denotes the subpath of $D_k \subseteq J$ with $V(X/R):= \{h(v)\,|\, v \in V(X)\}$ and 
$E(X/R):= \{h'(x)\,|\,x \in E(X) \cap (E(F)-E(R))\}$.

\begin{theorem}\label{pppp}
$H$ is a cyclically $5$-edge connected permutation snark.
\end{theorem}

Proof. By Theorem \ref{!!}, $F$ is not contained in a CDC. Hence Lemma \ref{bla} implies that $H$ is not $3$-edge colorable. 
It remains to show that  $\lambda_c (H) > 4$.

Suppose by contradiction that $E_0:=\{a_1,a_2,a_3,a_4\}$ is a cyclic $4$-edge cut of $H$ and let 
by Proposition \ref{perch} (2), w.l.o.g. $\{a_1,a_2\}$ be a matching of $E(C_1)$. Let $X$, $X'$ denote the two components of $C_1-a_1-a_2$.  $X/R$ and $X/R'$ are edge disjoint paths of $D_1 \subseteq J$ with $|E(X/R)| + |E(X'/R)| \leq 7$. Let w.l.o.g. $|E(X/R)| \leq 3$. Then, $2 \leq |V(X/R)| \leq 4$ since by Lemma \ref{kkj}, $\{a_1,a_2\} \not\subseteq E(A^i_j)$ for $i=1,2,3,4$ and $j=1,2,3$. Set $Y:= X/R$.

By Lemma \ref{4cc} (2), $X^N := \langle N(X) \cap V(C_2) \rangle $ is a path in $C_2$. Hence $Y^*:=X^N/R$ is a path of $D_2 \subseteq J$. Since only subpaths of $F \subseteq H$ are contracted in the construction of $J$ and since every vertex of $X \subseteq C_1$ is adjacent to a vertex of $X^N \subseteq C_2$ and vice versa, every vertex of $Y$ is adjacent to a vertex of $Y^*$ and vice versa. We make the following observation.

If $v$ is an endvertex of $Y$, then $|(N(v) \cap V(D_2)) \cap V(Y^*)| \geq 1$.\\
If $|V(Y)| > 2$ and $v$ is an inner vertex of $Y$, then $N(v) \cap V(D_2) \subseteq V(Y^*)$.

Set $U:= N(V(oY)) \cap V(D_2)$. Denote by $v_s$ and $v_t$ the two distinct endvertices of $Y$. Then, $Y^*$ = $\langle S \cup T \cup U\rangle  $ for some $S \subseteq N(v_s) \cap V(D_2)$ and for some $T \subseteq N(v_t) \cap V(D_2)$.

If $V(Y)=\{v_1,v_2\}$ then $V(Y^*) \in \{\{w_1,w_5\},\{w_1,w_3\},\{w_1,w_3,w_5\}\}$.
We abbreviate this conclusion by only writing the indices: $12 \rightsquigarrow 15,13,135$. We recall that $2 \leq |V(Y)| \leq 4$. Analogously, we consider for $Y$ the following cases where $Y \not= \,\langle \{v_6,v_7,v_1,v_2\} \rangle  $.

$23\rightsquigarrow 13,15,135$; $34 \rightsquigarrow 16$; $45 \rightsquigarrow 46$; $56 \rightsquigarrow 42,47,247$; $67 \rightsquigarrow 24,47,247$; $71 \rightsquigarrow 14$.

$123 \rightsquigarrow 135$; $234 \rightsquigarrow 136,156,356$; $345 \rightsquigarrow 146$; $456 \rightsquigarrow 2467,246,467$; $567 \rightsquigarrow 247$; 
$671 \rightsquigarrow 124,147,1247$; $712 \rightsquigarrow 134,145,1345$.

$1234 \rightsquigarrow 1356$; $2345 \rightsquigarrow 1346,1456,13456$; $3456 \rightsquigarrow 1246,1467,12467$; $4567 \rightsquigarrow 2467$;
$5671 \rightsquigarrow1247$; $7123 \rightsquigarrow 1345$.

In none of the above cases, $V(Y^*)$ can be the vertex set of a path in $D_2$ contradicting the assumption that $Y^*$ is a path.

Thus it remains to consider $Y:= \langle \{v_6,v_7,v_1,v_2 \} \rangle $. Since $v_1 \in V(Y) \subseteq V(D_1)$ (see the half-edges incident with $v_1$ in Figure 5), the path $X \subseteq C_1$ contains all vertices of $oA^1_1$ (Figure 1). The vertex $x^1_2 \in V(H)$ (Figure 1) corresponds to $v_3 \in V(D_1)$ (Figure 5) with $v_3 \not\in V(Y)$. Hence,  $x^1_2 \in V(C_1)$ and $x^1_2 \not\in V(X)$. Moreover, $x^1_2$ is neither adjacent to $x^1_4$ nor to $x^1_5$; otherwise $Q^1$ contains a circuit of length $4$ which contradicts Def. \ref{c5}. Hence, $\{x^1_4,x^1_5\} \subseteq V(X^N)$.

By the structure of $Q_1$ (Figure 1) and since $oA^1_1 \subseteq X$, $X^N$ contains every vertex of $oA^1_2-w$ where $w$ denotes the neighbor of $x^1_2$ in $oA^1_2$. Thus, and since $X^N$ is a subpath of $C_2$, and since $\{x^1_4,x^1_5\} \subseteq V(X^N)$, $|V(X^N)| = |V(C_2)|-1$. Since  $|V(X)| = |V(X^N)|$ and since $\{a_1,a_2\}$ in the definition of $X$ is a matching, this is impossible which finishes the proof.

Denote by $P_{10}$ the Petersen graph.

\begin{definition}\label{Hmenge}
 Set $\mathcal H:= \bigcup_{n=0}^\infty \{H_n\}$ where $H_n:=H(H_{n-1},P_{10},P_{10},P_{10})$ and 
$H_0:=P_{10}$, see Definition \ref{Hcon}.
\end{definition}

Note that in the above definition we use the graph $H_n$, respectively, $P_{10}$ as $Q^i$ and thus suppose that two subpaths in the known 
permutation $2$-factors of $H_n$ and $P_{10}$ are chosen as the paths specified in Definition \ref{c5}.

\begin{corollary}\label{wichtig}
$\mathcal H$ is an infinite set of cyclically $5$-edge connected permutation snarks where $H_n \in \mathcal H$ has $10+24n$ vertices.
\end{corollary}

Hence and by Theorem \ref{!!}, we obtain the following corollary.

\begin{corollary}
For every $n \in \mathbb N$, there is a counterexample of order $10+24n$ to Conjecture \ref{c2}, Conjecture \ref{c4} and Conjecture \ref{c5}.
\end{corollary}

\begin{corollary}
For every $n \in \mathbb N$, there is a counterexample $G$ of order $5+12n$ to Conjecture \ref{c1}.
\end{corollary}

Proof. Set $H:= H_n$, see Definition \ref{Hmenge}. Contract every spoke of $H$ with respect to $F$ to obtain a $4$-regular graph $G$ of order $5+12n$. Then $G= C'_1 \cup C'_2$ where $C'_1$ and $C'_2$ are two edge-disjoint hamiltonian circuits of $G$ which correspond to $C_1$ and $C_2$ in $H$.  Hence, every edge-cut $E_0$ of $G$ has even seize. Suppose that $E_0$ is an essential $4$-edge cut of $G$. Since $G$ is $4$-regular, every component of $G-E_0$ has more than $2$ vertices. It is straightforward to see that then $E_0$ corresponds to a cyclic $4$-edge cut of $H$ which contradicts $\lambda_c(H)=5$. Thus and since $E_0$ is of even size, $G$ is essentially $6$-edge connected. By defining that every pair of two edges which are adjacent and part of $C'_i$ for some $i \in \{1,2\}$ from a transition, we obtain a transition system $T(G)$ of $G$. Since every compatible cycle decomposition of $T(G)$ would imply a CDC $S$ of $H$ which contains $F$ and thus would contradict Theorem \ref{!!}, the proof is finished.

\textbf{Acknowledgement.}
A. Hoffmann-Ostenhof thanks R. H\"aggkvist for the invitation to the University of Umea where part of the work has been done.

\footnotesize

\end{document}